\documentclass[a4paper,12pt,final]{amsart}
\usepackage{times,a4wide,mathrsfs,amssymb,amsmath,amsthm,enumerate,xypic,tikzsymbols,dsfont}

\newcommand{\C}{\mathbb{C}}

\newcommand{\QQ}{\mathbb{Q}}
\newcommand{\NN}{\mathbb{N}}
\newcommand{\PP}{\mathbb{P}}

\newcommand{\OO}{\mathcal O}
\newcommand{\Ss}{\mathcal S}

\newcommand{\XX}{\mathcal X}
\newcommand{\YY}{\mathcal Y}

\newcommand{\MM}{\mathcal M}

\newcommand{\rom}{\romannumeral}

\newcommand{\one}{\mathds{1}}

\DeclareMathOperator{\ima}{Im}

\newtheorem{theorem}{Theorem}[section]
\newtheorem{claim}[theorem]{Claim}
\newtheorem{lemma}[theorem]{Lemma}

\newtheorem{corollary}[theorem]{Corollary}
\newtheorem{proposition}[theorem]{Proposition}

\newtheorem{convention}{Conventions}

\newtheorem{nonumbering}{Theorem}

\newtheorem{nonumberingc}{Corollary}

\theoremstyle{definition}
\newtheorem{remark}[theorem]{Remark}
\newtheorem{definition}[theorem]{Definition}

\newtheorem{notation}[theorem]{Notation}

\newtheorem{nonumberingt}{Acknowledgments}

\begin{document}

\author[Robert Laterveer]
{Robert Laterveer}

\address{Institut de Recherche Math\'ematique Avanc\'ee,
CNRS -- Universit\'e 
de Strasbourg,\
7 Rue Ren\'e Des\-car\-tes, 67084 Strasbourg CEDEX,
FRANCE.}
\email{robert.laterveer@math.unistra.fr}

\title[Algebraic cycles and intersections of three quadrics]{Algebraic cycles and intersections of three quadrics}

\begin{abstract} 
Let $Y$ be a smooth complete intersection of three quadrics, and assume the dimension of $Y$ is even. We show that $Y$ has a multiplicative Chow--K\"unneth decomposition, in the sense of Shen--Vial. 
As a consequence, the Chow ring of (powers of) $Y$ displays K3-like behaviour. As a by-product of the argument, we also establish a multiplicative Chow--K\"unneth decomposition for double planes.
\end{abstract}

\thanks{Supported by ANR grant ANR-20-CE40-0023.}


\maketitle

\section{Introduction}

Given a smooth projective variety $Y$ over $\C$, let $A^i(Y):=CH^i(Y)_{\QQ}$ denote the Chow groups of $Y$ (i.e. the groups of codimension $i$ algebraic cycles on $Y$ with $\QQ$-coefficients, modulo rational equivalence). The intersection product defines a ring structure on $A^\ast(Y)=\bigoplus_i A^i(Y)$, the {\em Chow ring\/} of $Y$ \cite{F}.
In the case of K3 surfaces, this ring structure has a peculiar property:

\begin{theorem}[Beauville--Voisin \cite{BV}]\label{bv} Let $S$ be a projective K3 surface. Then
   \[ \ima\Bigl( A^1(S)\otimes A^1(S)\to A^2(S)\Bigr) = \QQ[c_2(S)]\ .\]
   \end{theorem}

Motivated by Theorem \ref{bv}, Beauville \cite{Beau3} has conjectured that for certain special varieties, the Chow ring should admit a multiplicative splitting. To make concrete sense of Beauville's elusive ``splitting property'' conjecture, Shen--Vial \cite{SV} introduced the concept of {\em multiplicative Chow--K\"unneth decomposition\/}. It is something of a challenge to understand the class of special varieties admitting such a decomposition: abelian varieties, K3 surfaces and cubic hypersurfaces are in this class, hyperelliptic curves are in this class (but not all curves), some Fano varieties are in this class (but not all Fano varieties), cf. subsection \ref{ss:mck} below.

The aim of the present paper is to find new members of this class, by considering complete intersections of quadrics. We note that the case of one quadric is trivial (quadrics have trivial Chow groups), while the case of intersections of two quadrics was treated in \cite{2quadrics}. The main result of this paper is about three quadrics:

\begin{nonumbering}[=Theorem \ref{main} and Corollary \ref{cor}] Let $Y\subset\PP^{n+3}(\C)$ be a smooth dimensionally transverse intersection of three quadrics, where $n\ge 2$ is even. Then $Y$ has a multiplicative Chow--K\"unneth decomposition. In particular, 
     \[  \ima\Bigl( A^i(Y)\otimes A^j(Y)\to A^{i+j}(Y)\Bigr) =\QQ[c_{i+j}(Y)]\ \ \ \ \ \forall\  i,j>0\ \]
     \end{nonumbering}

In proving Theorem \ref{main}, we rely on the connection between $Y$ and a double plane \cite{OG}, and we apply instances of the {\em Franchetta property\/} (cf. Section \ref{s:fr} below). 

%
%

Using Theorem \ref{main}, we also prove a result concerning the {\em tautological ring\/}, which is a certain subring of the Chow ring of powers of $Y$:

\begin{nonumberingc}[=Corollary \ref{cor1}] Let $Y\subset\PP^{n+3}(\C)$ be as in Theorem \ref{main}, and $m\in\NN$. Let
  \[ R^\ast(Y^m):=\Bigl\langle (p_i)^\ast(h), (p_{ij})^\ast(\Delta_Y)\Bigr\rangle\ \subset\ \ \ A^\ast(Y^m)   \]
  be the $\QQ$-subalgebra generated by pullbacks of the polarization $h\in A^1(Y)$ and pullbacks of the diagonal $\Delta_Y\in A^{n}(Y\times Y)$. 
  Then the cycle class map induces injections
   \[ R^\ast(Y^m)\ \hookrightarrow\ H^\ast(Y^m,\QQ)\ \ \ \ \forall\ m\le 2\dim H^{n}(Y,\QQ)-1  \ .\]
   
 Moreover, $R^\ast(Y^m)\to H^\ast(Y^m,\QQ)$ is injective for all $m$ if and only if $Y$ is Kimura finite-dimensional \cite{Kim}.
     \end{nonumberingc}

This is similar to existing results for hyperelliptic curves and for K3 surfaces (cf. \cite{Ta2}, \cite{Ta}, \cite{Yin} and Remark \ref{tava} below).

 Theorem \ref{main} also has consequences for the Chow ring of double planes:
 
 \begin{nonumbering}[=Theorem \ref{main2}] Let $S$ be the double cover of $\PP^2$ branched along a smooth curve $C\subset\PP^2$. Then $S$ has a multiplicative Chow--K\"unneth decomposition.
 In particular,
   \[ \ima\Bigl( A^1(S)\otimes A^1(S)\to A^2(S)\Bigr) = \QQ[c_2(S)]\ .\]
   \end{nonumbering}
   
 Let us end this introduction with some open questions: do complete intersections of three quadrics in an {\em even-dimensional\/} projective space have a multiplicative Chow--K\"unneth decomposition ? Do complete intersections of {\em four or more\/} quadrics have a multiplicative Chow--K\"unneth decomposition ? The argument proving Theorem \ref{main} does not allow to settle these questions (cf. Remark \ref{difficult} below).

 \vskip0.6cm

\begin{convention} In this article, the word {\sl variety\/} will refer to a reduced irreducible scheme of finite type over $\C$. A {\sl subvariety\/} is a (possibly reducible) reduced subscheme which is equidimensional. 

{\bf All Chow groups are with rational coefficients}: we will denote by $A_j(Y)$ the Chow group of $j$-dimensional cycles on $Y$ with $\QQ$-coefficients; for $Y$ smooth of dimension $n$ the notations $A_j(Y)$ and $A^{n-j}(Y)$ are used interchangeably. 
The notation $A^j_{hom}(Y)$ will be used to indicate the subgroup of homologically trivial cycles.
For a morphism $f\colon X\to Y$, we will write $\Gamma_f\in A_\ast(X\times Y)$ for the graph of $f$.

The contravariant category of Chow motives (i.e., pure motives with respect to rational equivalence as in \cite{Sc}, \cite{MNP}) will be denoted 
$\MM_{\rm rat}$.
\end{convention}

 \vskip0.6cm

\goodbreak
\section{Preliminaries}

\subsection{Intersections of three quadrics} 

\begin{theorem}[O'Grady \cite{OG}]\label{net} Let $\Lambda\cong\PP^2$ be a net of quadrics in $\PP^{n+3}$ where $n\ge 2$ is even, and such that the base locus of the net $Y\subset\PP^{n+3}$ is smooth. Assume the degree $n+4$ curve $C\subset\Lambda\cong\PP^2$
parametrizing singular quadrics in the net is smooth, and let $\pi\colon S\to\PP^2$ be the double cover branched along $C$.

\noindent
(\rom1) There is a correspondence $\Gamma\in A^{{n\over 2}+1}(Y\times S)$ inducing an isomorphism
  \[ \Gamma_\ast\colon\ \ H^n_{prim}(Y,\QQ)\ \xrightarrow{\cong}\ H^2_{prim}(S,\QQ)\ .\]
  (Here $H^n_{prim}(Y,\QQ)$ is the kernel of the map $H^n(Y,\QQ)\to H^{n+2}(Y,\QQ)$ given by cupping with an ample class, and
  $H^2_{prim}(S,\QQ)$ is defined as the anti-invariant part with respect to the covering involution of $S$.)

  \noindent
  (\rom2) The map $\Gamma_\ast$ of (\rom1) is compatible with cup-product:
    \[    \Gamma_\ast(\alpha)\cup \Gamma_\ast(\beta) = \alpha\cup\beta\ \ \ \forall\ \alpha,\beta\ \in\ H^n_{prim}(Y,\QQ)\ .\]
\end{theorem}

\begin{proof} This is the main result of \cite{OG}, extending earlier work of Mukai for the case $n=2$ (where $Y$ and $S$ are both K3 surfaces, related as moduli spaces) \cite{Muk}. Part (\rom1) is also proven (and generalized to complete intersections of an arbitrary number of quadrics in an odd-dimensional projective space) by Terasoma \cite[Corollary 4.5.3]{Tera}.

Let us briefly explain the correspondence $\Gamma$ used by O'Grady. We consider the incidence variety
  \[ X:=\Bigl\{ (\lambda,p)\in \Lambda\times\PP^{n+3}\, \Big\vert\, p\in Q_\lambda\Bigr\}\ \ \ \subset\ \Lambda\times\PP^{n+3}\ ,\]
  of dimension $\dim X=n+4$ (this is the total space of the net $\Lambda$).
Let us write $\iota\colon \Lambda\times Y\hookrightarrow X$ for the inclusion morphism, and $p_2\colon \Lambda\times Y\to Y$ for the projection. It is readily proven that
  \begin{equation}\label{one}  \iota_\ast (p_2)^\ast\colon\ \ H^n_{prim}(Y,\QQ)\ \to\ H^{n+4}_{prim}(X,\QQ) \end{equation}
  is an isomorphism which is compatible with cup-product (this is \cite[Theorem 1.1]{OG}).
  
 To link $X$ to the double plane $S$, one introduces the variety of ``relative linear subspaces''
   \[ P:= \Bigl\{  M\subset X\, \Big\vert\, (p_1)(M)\ \hbox{is\ a\ point\ and\ $(p_2)(M)\subset\PP^{n+3}$\ is\ a\ $\PP^{{n\over 2}+1}$}\Bigr\}\ ,\]
   and the incidence variety
   \[ I_P:= \Bigl\{ (p,x)\in P\times X\,\Big\vert\,  x\in M_p\Bigr\}\ \ \ \subset\ P\times X\ ,\]
   where $M_p\subset X$ is the subvariety corresponding to $p\in P$. 
  This variety $P$ has the property that the Stein factorization of the natural morphism $P\to\Lambda$ is 
    \[ P\ \xrightarrow{f}\  S\ \xrightarrow{\pi}\ \Lambda\ .\]
   
Let $d:=2 \dim P$. The correspondence $I_P$ induces a map
 \[ (I_P)_\ast\colon\ \ H^{d-2}(P,\QQ)\ \to\ H^{n+4}(X,\QQ) \]
 (this is the map labelled $\theta$ in \cite[Definition 2.1]{OG}; it is an isomorphism when restricted to the primitive parts). It is proven in loc. cit. that the composition
  \begin{equation}\label{two} H^2_{prim}(S,\QQ)\ \xrightarrow{f^\ast}\ H^2_{prim}(P,\QQ)\ \xrightarrow{\cdot h_P^{{d\over 2}-2}}\ H^{d-2}_{prim}(P,\QQ)\ \xrightarrow{(I_P)_\ast}\ H^{n+4}_{prim}(X,\QQ) \end{equation}
  (where $h_P\in A^1(P)$ denotes an ample class)
  is an isomorphism, and compatible with cup-product.
  
  Combining \eqref{one} and \eqref{two}, one obtains two correspondences $\Gamma_1$ and $\Gamma_2$ inducing isomorphisms
   \[ \begin{array}[c]{ccccc} H^n_{prim}(Y,\QQ)& \xrightarrow{(\Gamma_1)_\ast}&  H^{n+4}_{prim}(X,\QQ)& \xleftarrow{(\Gamma_2)_\ast}& H^2_{prim}(S,\QQ)\ .\\
                                                                             &\cong    &&\cong&\\
                                                                             \end{array}\]
 The fact that $(\Gamma_2)_\ast$ is compatible with cup-product implies that the inverse to $(\Gamma_2)_\ast$ is induced by the transpose correspondence ${}^t \Gamma_2$  
 (cf. for instance \cite[Lemma 5]{V43}).
We obtain the required correspondence by composition: 
   \[ \Gamma:=  {}^t \Gamma_2\circ\Gamma_1\ \ \ \in\ A^{{n\over 2}+1}(Y\times S)\ .\]
   \end{proof}
 
 \begin{remark} The cohomological relation of Theorem \ref{net} also exists (and in fact, is explained by) a relation on the level of derived categories \cite{BO}, \cite{Ku0}, \cite{Add}: given $Y$ and $S$ as in Theorem \ref{net},
 there is a semi-orthogonal decomposition
    \[ D^b(Y)=\Bigl\langle \OO_Y(-n+3),\ldots,\OO_Y(-1), \OO_Y, D^b(S,\alpha)\Bigr\rangle\ ,\]
    where $\alpha\in H^2(S,\OO^\ast_S)$ is some Brauer class.
    
 In the present paper, we will {\em not\/} use this categorical relation; the main input for our argument will be the cohomological relation of Theorem \ref{net}.   
 \end{remark}

 \subsection{MCK decomposition}
\label{ss:mck}

\begin{definition}[Murre \cite{Mur}] Let $X$ be a smooth projective variety of dimension $n$. We say that $X$ has a {\em CK decomposition\/} if there exists a decomposition of the diagonal
   \[ \Delta_X= \pi^0_X+ \pi^1_X+\cdots +\pi_X^{2n}\ \ \ \hbox{in}\ A^n(X\times X)\ ,\]
  such that the $\pi^i_X$ are mutually orthogonal idempotents and $(\pi_X^i)_\ast H^\ast(X,\QQ)= H^i(X,\QQ)$.
  
  (NB: ``CK decomposition'' is shorthand for ``Chow--K\"unneth decomposition''.)
\end{definition}

\begin{remark} The existence of a CK decomposition for any smooth projective variety is part of Murre's conjectures \cite{Mur}, \cite{J4}. 
\end{remark}

\begin{definition}[Shen--Vial \cite{SV}] Let $X$ be a smooth projective variety of dimension $n$. Let $\Delta_X^{sm}\in A^{2n}(X\times X\times X)$ be the class of the small diagonal
  \[ \Delta_X^{sm}:=\bigl\{ (x,x,x)\ \vert\ x\in X\bigr\}\ \subset\ X\times X\times X\ .\]
  An {\em MCK decomposition\/} is a CK decomposition $\{\pi_X^i\}$ of $X$ that is {\em multiplicative\/}, i.e. it satisfies
  \[ \pi_X^k\circ \Delta_X^{sm}\circ (\pi_X^i\times \pi_X^j)=0\ \ \ \hbox{in}\ A^{2n}(X\times X\times X)\ \ \ \hbox{for\ all\ }i+j\not=k\ .\]
  
 (NB: ``MCK decomposition'' is shorthand for ``multiplicative Chow--K\"unneth decomposition''.) 
  
  \end{definition}
  
  \begin{remark} The small diagonal (seen as a correspondence from $X\times X$ to $X$) induces the {\em multiplication morphism\/}
    \[ \Delta_X^{sm}\colon\ \  h(X)\otimes h(X)\ \to\ h(X)\ \ \ \hbox{in}\ \MM_{\rm rat}\ .\]
 Suppose $X$ has a CK decomposition
  \[ h(X)=\bigoplus_{i=0}^{2n} h^i(X)\ \ \ \hbox{in}\ \MM_{\rm rat}\ .\]
  By definition, this decomposition is multiplicative if for any $i,j$ the composition
  \[ h^i(X)\otimes h^j(X)\ \to\ h(X)\otimes h(X)\ \xrightarrow{\Delta_X^{sm}}\ h(X)\ \ \ \hbox{in}\ \MM_{\rm rat}\]
  factors through $h^{i+j}(X)$.
  
  If $X$ has an MCK decomposition, then setting
    \[ A^i_{(j)}(X):= (\pi_X^{2i-j})_\ast A^i(X) \ ,\]
    one obtains a bigraded ring structure on the Chow ring: that is, the intersection product sends $A^i_{(j)}(X)\otimes A^{i^\prime}_{(j^\prime)}(X) $ to  $A^{i+i^\prime}_{(j+j^\prime)}(X)$.
    
      It is expected that for any $X$ with an MCK decomposition, one has
    \[ A^i_{(j)}(X)\stackrel{??}{=}0\ \ \ \hbox{for}\ j<0\ ,\ \ \ A^i_{(0)}(X)\cap A^i_{hom}(X)\stackrel{??}{=}0\ ;\]
    this is related to Murre's conjectures B and D, that have been formulated for any CK decomposition \cite{Mur}.

  The property of having an MCK decomposition is severely restrictive, and is closely related to Beauville's ``splitting property conjecture'' \cite{Beau3}. 
  To give an idea: hyperelliptic curves have an MCK decomposition \cite[Example 8.16]{SV}, but the very general curve of genus $\ge 3$ does not have an MCK decomposition \cite[Example 2.3]{FLV2}. As for surfaces: a smooth quartic in $\PP^3$ has an MCK decomposition, but a very general surface of degree $ \ge 7$ in $\PP^3$ should not have an MCK decomposition \cite[Proposition 3.4]{FLV2}.
For more detailed discussion, and examples of varieties with an MCK decomposition, we refer to \cite[Section 8]{SV}, as well as \cite{V6}, \cite{SV2}, \cite{FTV}, \cite{37}, \cite{38}, \cite{39}, \cite{40}, \cite{FLV2}, \cite{44}, \cite{g8}, \cite{2quadrics}.
   \end{remark}

 \section{The Franchetta property}
 \label{s:fr}
 
 \subsection{Definition}
 
 \begin{definition} Let $\XX\to B$ be a smooth projective morphism, where $\XX, B$ are smooth quasi-projective varieties. We say that $\XX\to B$ has the {\em Franchetta property in codimension $j$\/} if the following holds: for every $\Gamma\in A^j(\XX)$ such that the restriction $\Gamma\vert_{X_b}$ is homologically trivial for the very general $b\in B$, the restriction $\Gamma\vert_{X_b}$ is zero in $A^j(X_b)$ for all $b\in B$.
 
 We say that $\XX\to B$ has the {\em Franchetta property\/} if $\XX\to B$ has the Franchetta property in codimension $j$ for all $j$.
 \end{definition}
 
 This property is studied in \cite{PSY}, \cite{BL}, \cite{FLV}, \cite{FLV3}.
 
 \begin{definition} Given a family $\XX\to B$ as above, with $X:=X_b$ a fiber, we write
   \[ GDA^j_B(X):=\ima\Bigl( A^j(\XX)\to A^j(X)\Bigr) \]
   for the subgroup of {\em generically defined cycles}. 
   (In a context where it is clear to which family we are referring, the index $B$ will sometimes be suppressed from the notation.)
  \end{definition}
  
  With this notation, the Franchetta property amounts to saying that $GDA^\ast(X)$ injects into cohomology, under the cycle class map. 
 
\subsection{The families} We introduce various families that will be used in the argument.

  \begin{notation}\label{not} Given $n\in\NN$ even and at least 2, let
   \[ B\ \subset\ \bar{B}:=\PP H^0\bigl(\PP^{n+3},\OO_{\PP^{n+3}}(2) \oplus  \OO_{\PP^{n+3}}(2) \oplus\OO_{\PP^{n+3}}(2)\bigr)\] 
   denote the Zariski open parametrizing smooth dimensionally transverse complete intersections of three quadrics, and let 
   \[ \YY\ \to\ B \]
   denote the universal family of smooth $n$-dimensional complete intersections of three quadrics. 
   
   Let
     \[ \MM_{(2,2,2)}:= B/ PGL(n+3) \]
     be the moduli space, and $\YY\to\MM_{(2,2,2)}$ the universal complete intersection of three quadrics (these exist as algebraic stacks).
  \end{notation}
  
  \begin{notation}\label{not2} Given $n\in\NN$ even and at least 2, let us write $P:= \PP(1,1,1,{n\over 2}+2)$ for the weighted projective space. Let   
    \[  B_{dp}\ \subset\ \bar{B}_{dp}:= \PP H^0\bigl(P^{},\OO_{P^{}}(n+4)\bigr)\] 
   denote the Zariski open parametrizing smooth surfaces, and let 
   \[ \Ss\ \to\ B_{dp} \]
   denote the universal family of smooth hypersurfaces. 
    (NB: $\Ss \to B_{dp}$ may be identified with the universal family of double planes branched over curves of degree $n+4$; the index ``dp'' is shorthand for ``double planes''.)
    
    Let $\MM_{dp}$ be the moduli space of double planes branched over degree $n+4$ curves, and $\Ss\to\MM_{dp}$ the universal double plane (in the sense of stacks).
     \end{notation}
   
  \begin{lemma}\label{F} The universal families $\YY$ (lying over $B$) and $\Ss$ (lying over $B_{dp}$) are smooth.
   \end{lemma}
   
   \begin{proof}   Let $\bar{\YY}\supset \YY$ denote the Zariski closure of $\YY$ in $\PP^n\times \bar{B}$.
  As $\OO_{\PP^{n}}(3)\oplus \OO_{\PP^n}(2)$ is base-point free, each point of $\PP^n$ imposes 2 independent conditions on $\bar{B}$. This means that the projection $\YY\to\PP^n$ has the structure of a projective bundle, hence is smooth.  
   
 The argument for $\Ss$ is similar.
     \end{proof}

  The following result relates (the parameter spaces $B$ and $B_{dp}$ and) the moduli spaces $\MM_{(2,2,2)}$ and $\MM_{dp}$:
 
 \begin{proposition}\label{key} There are morphisms
   \[   \begin{array}[c]{ccc}
         B && B_{dp}\\
         &&\\
         \downarrow&&\downarrow\\
         &&\\
         \MM_{(2,2,2)} & \stackrel{\phi}{\dashrightarrow}& \ \MM_{dp}\ .\\
         \end{array}\]
         The vertical arrows are surjections. The arrow $\phi$ is surjective, and there exists a non-empty Zariski open $\MM^0_{(2,2,2)}\subset\MM_{(2,2,2)}$ such that $\phi\colon \MM^0_{(2,2,2)}\to \MM_{dp}$ is a bijection.     
         \end{proposition}
         
         \begin{proof} The map $\phi$ arises from Theorem \ref{net}; it sends a net of quadrics to the double plane branched along the Hesse curve $C$ (associated to the net of quadrics).
         The map $\phi$ is well-defined on the locus where the Hesse curve $C$ is non-singular.
         The fact that $\phi$ is surjective is proven in \cite[Proposition 4.2(a) and Remark 4.4]{Bdet} (according to loc. cit., this is well-known and should be attributed to Dixon \cite{Dix}). The fact that $\phi$ is generically injective follows from the generic Torelli theorem for intersections of three quadrics in odd-dimensional projective space \cite[Th\'eor\`eme IV.4.1]{Las}.
           \end{proof}

\begin{remark} Proposition \ref{key} says that every double plane arises as a double plane associated to a net of quadrics. This fact will be essential for proving the main result of the present paper (Theorem \ref{main}). The same fact is also used in proving generic Torelli for intersections of three quadrics in odd-dimensional projective space \cite{Las}, by reducing to generic Torelli for double planes.
\end{remark}

  \subsection{Franchetta for $Y$}
  
 \begin{proposition}\label{f1} Let $\YY\to B$ be the universal family of smooth complete intersections of multidegree $(2,2,2)$ in $\PP^{n+3}$ (Notation \ref{not}).
 The family $\YY\to B$ has the Franchetta property.
 \end{proposition}
 
 \begin{proof} This is easy, and does not need the assumption $n$ even. We have already seen (proof of Lemma \ref{F}) that $\bar{\YY}\to\PP^{n+3}$ is a projective bundle. Using the projective bundle formula, and reasoning as in
 \cite{PSY} and \cite{FLV}, this implies that
   \[ GDA^\ast_B(Y)=\ima\bigl( A^\ast(\PP^{n+3})\to A^\ast(Y)\bigr)=\langle h\rangle\ ,\]
   where $h\in A^1(Y)$ is the restriction of a hyperplane section. It follows that $GDA^\ast_B(Y)$ injects into cohomology.
   \end{proof}

  \subsection{Franchetta for $Y^2$}
  
  \begin{notation} Given a morphism $\XX\to B$, let us write 
    \[ \XX^{n/B}:=\XX\times_B \XX\times_B \cdots \times_B \XX \]
    for the $n$th relative power.
   \end{notation}

   \begin{proposition}\label{f2mid} Let $\YY\to B$ be as above. Then the family $\YY^{2/B}\to B$ has the Franchetta property in codimension $n$.
  \end{proposition}  
  
  \begin{proof} We note that $\OO_{\PP^{n+3}}(2)$ is (not only base-point free but even) very ample, which means that the set-up verifies condition $(\ast_2)$ of \cite{FLV}. Then applying the stratified projective bundle argument as in loc. cit. (in the precise form of \cite[Proposition 2.6]{FLV3}), we find that
    \[ GDA^\ast(Y^2)=\langle h, \Delta_Y\rangle\ .\]
    Let us now check that $GDA^n(Y^2)$ injects into cohomology. One observes that the class of the diagonal in $H^{2n}(Y\times Y,\QQ)$ can not be expressed as a decomposable cycle (for otherwise $Y$ would not have any transcendental cohomology, which is absurd), and so the required injectivity follows from Proposition \ref{f1}.  
  \end{proof}
  
  On a similar note, we have:
  
  \begin{proposition}\label{s2mid} Let $\Ss\to B_{dp}$ be as above. Then the family $\Ss^{2/B_{dp}}\to B_{dp}$ has the Franchetta property in codimension $\le 3$.
  \end{proposition}  
  
  \begin{proof} As above, we write $\PP:=\PP(1,1,1,{n\over 2}+2)$ for the weighted projective space. We know that the line bundle $\OO_\PP(n+4)$ is very ample on $\PP$ (this follows from \cite[Proposition 2.3(\rom3)]{Del}). Once more applying \cite[Proposition 2.6]{FLV3}, it follows that
    \begin{equation}\label{know}GDA^\ast_{B_{dp}}(S\times S)=\langle h,\Delta_S\rangle\ ,\end{equation}
   where $h\in A^1(S)$ is the restriction of a hyperplane section. Since $S$ has non-zero transcendental cohomology, the injectivity of $GDA^2(S\times S)$ follows as above.
   The excess intersection formula implies that there is an equality
       \begin{equation}\label{excs} \Delta_S\cdot (p_j)^\ast(h) =\sum_j  a_j(p_1)^\ast(h^j)\cdot (p_2)^\ast(h^{3-j})\ \ \ \hbox{in}\ A^3(S\times S)\ ,\end{equation} 
   and so the injectivity of
   $GDA^3(S\times S)\to H^6(S\times S,\QQ)$ follows from the injectivity of $\langle h\rangle\to H^\ast(S,\QQ)$ (which is elementary).  
    \end{proof}

\begin{definition} Given any $n$-dimensional smooth complete intersection $X$ in a weighted projective space, one can define the {\em primitive motive\/} $h^n_{prim}(X):=(X,\pi^{n,prim}_X,0)\in\MM_{\rm rat}$, where
   \[ \pi^{n,prim}_X:= \Delta_X - {1\over d}\, \sum_{j=0}^n  h^j\times h^{n-j}\ \ \ \in\ A^n(X\times X)\ .\]
   (Here $h\in A^1(X)$ is a hyperplane section, and $d$ is the degree of $X$.)
  \end{definition}

  As a consequence of Propositions \ref{f2mid} and \ref{s2mid}, we can upgrade O'Grady's result (Theorem \ref{net}) to the level of Chow motives:
  
  \begin{corollary}\label{netchow} Let $Y\subset\PP^{n+3}$ and $S$ and $\Gamma$ be as in Theorem \ref{net}. There is an isomorphism of Chow motives
     \[ \Gamma\colon\ \ h^n_{prim}(Y)\ \xrightarrow{\cong}\ h^2_{prim}(S)(1-{n\over 2})   \ \ \ \hbox{in}\ \MM_{\rm rat}\ .\]
       \end{corollary}
       
  \begin{proof} Looking at the proof of Theorem \ref{net}, one sees that the correspondence $\Gamma$ is generically defined (with respect to $B$). Because $\Gamma$ is compatible with cup-product, the inverse to the isomorphism $\Gamma_\ast$ of Theorem \ref{net} is induced by the transpose ${}^t \Gamma$. That is, one has an equality of correspondences in cohomology
  \[ \pi^{n,prim}_Y =    \pi^{n,prim}_Y\circ {}^t \Gamma \circ \pi^{2,prim}_S\circ \Gamma\circ    \pi^{n,prim}_Y\ \ \ \hbox{in}\ H^{2n}(Y\times Y,\QQ)\ .\]
  Applying Proposition \ref{f2mid}, the same equality already holds in $A^n(Y\times Y)$, and so there is a split injection
    \[ \Gamma\colon\ \  h^n_{prim}(Y)\ \hookrightarrow\ h^2_{prim}(S)\ \ \ \hbox{in}\ \MM_{\rm rat}\ .\]
    
    Next, we note that the construction of the correspondence $\Gamma$ can also be done on the level of the moduli spaces $\MM_{(2,2,2)}$ and $\MM_{dp}$.
    As such the cycle
    \[  R:= \pi^{2,prim}_S - \pi^{2,prim}_S\circ  {} \Gamma \circ \pi^{n,prim}_Y\circ {}^t \Gamma\circ \pi^{2,prim}_S\ \ \ \in\ A^2(S\times S) \]
    is generically defined, with respect to $\MM_{dp}$. A fortiori, $R$ is generically defined with respect to the parameter space $B_{dp}$. The cycle $R$ is homologically trivial (Theorem \ref{net}), and so thanks to Proposition \ref{s2mid} we find that $R=0$ in $A^2(S\times S)$. This proves the corollary.
      \end{proof}

  \begin{proposition}\label{f2} Let $\YY\to B$ be as above. Then the family $\YY^{2/B}\to B$ has the Franchetta property.
  \end{proposition}
  
  \begin{proof} We have already seen (in the proof of Proposition \ref{f2mid}) that
      \[ GDA^\ast(Y^2)=\langle h, \Delta_Y\rangle\ ,\]
      and (in view of Proposition \ref{f2mid}) we only need to check that $GDA^j(Y^2)$ injects into cohomology, for $j\not=n$.
    Let us first dispose of the case $j=2n$. In case $n=2$, $Y$ is a K3 surface and so $\Delta_Y^2$ can be expressed as a multiple of $(p_1)^\ast(h^2)\cdot (p_2)^\ast(h^2)$ \cite{BV}.
    In case $n>2$, $Y$ is Fano and so we do not need to worry about $A^{2n}(Y^2)=\QQ[h^{2n}]$.
  Using Lemma \ref{ok} below, $GDA^j(Y^2)$ for $n<j<2n$ is {\em decomposable\/}, i.e. 
    \[ GDA^j(Y^2)\ \ \subset\ GDA^\ast(Y)\otimes GDA^\ast(Y)\ \ \ \forall\ n< j<2n\ .\]
    Thus for $n< j<2n$, the injectivity reduces to Proposition \ref{f1}.
        
    \begin{lemma}\label{ok} Let $Y\subset\PP^{n+3}$ be a smooth complete intersection of three quadrics, where $n$ is even. Then
    \[  \Delta_Y\cdot (p_i)^\ast(h) = {\displaystyle\sum_k}\ {1\over 8}\, (p_1)^\ast (h^k)\cdot (p_2)^\ast (h^{n+1-k})\ \ \ \hbox{in}\ A^{n+1}(Y\times Y)\ .\]
    \end{lemma}
    
    To prove the lemma, we start by observing that the equality of the lemma holds true in cohomology: indeed, both sides of the equality act on $H^\ast(Y,\QQ)$ as multiplication with $h$ (note that the degree of $Y$ is $8$), and the same remains true when going to $Y\times M$ (for an arbitrary smooth projective variety $M$) and looking at the action on $H^\ast(Y\times M,\QQ)$.
    Next, we note that all cycles occurring in the lemma are generically defined, not only with respect to the parameter space $B$ but also with respect to the moduli space $\MM_{(2,2,2)}$.
    That is, to prove the lemma we need that $GDA^{n+1}_{\MM_{(2,2,2)}}(Y^2)$ injects into cohomology. We may assume $Y=Y_b$ for $b\in \MM^0_{(2,2,2)}$. Using Proposition \ref{key},
    we are reduced to proving that $GDA^{n+1}_{\MM_{dp}}(Y^2)$ injects into cohomology. Using Corollary \ref{netchow}, this is equivalent to the injectivity of
      \[  GDA^3_{\MM_{dp}}(S^2)\ \to\ H^6(S^2,\QQ) \ ,\]
      where $S$ is the double plane associated to $Y$. This follows from the injectivity of
       \[  GDA^3_{B_{dp}}(S^2)\ \to\ H^6(S^2,\QQ) \ ,\]  
       which is Proposition \ref{s2mid}.
           
      The lemma, and the proposition, are now proven. 
   \end{proof}
   
   \begin{remark}\label{!} The equality of Lemma \ref{ok} is remarkable, because this equality may fail for other complete intersections.
  Indeed, in case $Y\subset\PP^{m}$ is a smooth {\em hypersurface\/} (of any degree), there is equality
    \begin{equation}\label{form} \Delta_Y\cdot (p_i)^\ast(h)= {\displaystyle\sum_k} \ a_k\, (p_1)^\ast (h^k)\cdot (p_2)^\ast (h^{m-k})\ \ \ \hbox{in}\ A^{m}(Y\times Y)\ ,\end{equation}
    with $a_k\in\QQ$,
as follows from the excess intersection formula. On the other hand, in case $Y\subset\PP^m$ is a complete intersection of codimension at least 2,
in general there is {\em no equality\/} of the form \eqref{form}. Indeed, let $C$ be a very general curve of genus $g\ge 4$. The Faber--Pandharipande cycle
  \[  FP(C):= \Delta_C\cdot (p_j)^\ast(K_C) - {1\over 2g-2} K_C\times K_C\ \ \ \in A^2(C\times C)\ \ \ \ \ (j=1,2) \]
  is homologically trivial but non-zero in $A^2(C\times C)$ \cite{GG}, \cite{Yin0} (this cycle $FP(C)$ is the ``interesting 0-cycle'' in the title of \cite{GG}). In particular, for the very general complete intersection $Y\subset\PP^3$ of bidegree $(2,3)$, the cycle
  \[ FP(Y):=  \Delta_Y\cdot (p_j)^\ast(h) - {1\over 6} h \times h\ \ \ \in A^2(Y\times Y) \]
  is homologically trivial but non-zero, and so there cannot exist an equality of the form \eqref{form} for $Y$.
  
  This is intimately related to MCK decompositions. Indeed, if the curve $Y$ has an MCK decomposition which is generically defined, then $FP(Y)\in A^2_{(0)}(Y\times Y)\cong\QQ$ and so
  $FP(Y)$ would be zero.
  \end{remark}

 \begin{proposition}\label{s2} The family $\Ss^{2/B_{dp}}\to B_{dp}$ has the Franchetta property.  
  \end{proposition}
  
  \begin{proof} In view of Proposition \ref{s2mid}, it only remains to prove that $GDA^4_{B_{dp}}(S^2)$ (is 1-dimensional and so) injects into cohomology.    
  In view of \eqref{know}, it suffices to prove that
     \[ \Delta_S^2 =\lambda (p_1)^\ast(h^2)\cdot (p_2)^\ast(h^2)\ \ \ \hbox{in}\ A^4(S\times S)\ ,\]
     for some $\lambda\in\QQ$. That is, we are reduced to proving that $GDA^4_{\MM_{dp}}(S^2)$ is 1-dimensional.
     In view of Proposition \ref{key} and Corollary \ref{netchow}, we have
       \[  GDA^4_{\MM_{dp}}(S^2)\ \hookrightarrow\  GDA^{n+2}_{\MM_{dp}}(Y^2)\ \cong\ GDA^{n+2}_{\MM^0_{(2,2,2)}}(Y^2)\ .\]
    The result now follows from the Franchetta property for $\YY^{2/B}\to B$ (Proposition \ref{f2}).
  \end{proof}

  \subsection{Franchetta for $Y^3$ in codimension $\le 2n$}
  
  \begin{proposition}\label{s3} Let $\Ss\to B_{dp}$ be as in Notation \ref{not2}. The family $\Ss^{3/B_{dp}}\to B_{dp}$ has the Franchetta property in codimension $\ge 5$.
      \end{proposition}
      
    \begin{proof} We observe that the set-up $(\PP,\OO_\PP({n\over 2}+2))$ has the property $(\ast_3)$ of \cite{FLV}: for $n=2$ this was verified in \cite{FLV}, the general case follows by the same argument (3 distinct points on a line impose 3 independent conditions on the linear system). The argument of \cite{FLV} (in the precise form of \cite[Proposition 2.6]{FLV3}) then gives an equality
       \[  GDA^\ast_{B_{dp}}(S^3)=\langle  h,\Delta_S\rangle\ .\]
       
 For the case of codimension 5, we remark that an element in $GDA^5(S^3)$ is thus a linear combination of elements that are decomposable (in the sense that they come from $GDA^\ast(S^2)\otimes GDA^\ast(S)$, and elements of the form
      \[  \Delta_S^{sm}\cdot(p_j)^\ast(h)=  (p_{12})^\ast(\Delta_S)\cdot (p_{23})^\ast(\Delta_S)\cdot (p_j)^\ast(h)\ .\]
    Using \eqref{excs}, such an element can be expressed in terms of one diagonal and the $(p_i)^\ast(h)$ and so we find that
     \[ GDA^5(S^3)\ \subset\ \Bigl( GDA^\ast(S^2)\otimes GDA^\ast(S)\Bigr)^{\oplus 3}\ .\]
     Using the K\"unneth decomposition in cohomology, plus the injectivity of $GDA^\ast(S^2)\to H^\ast(S^2,\QQ)$ (Proposition \ref{s2}), we find that $GDA^5(S^3)$ injects into cohomology.
     
   For the case of codimension 6, we remark that an element in $GDA^6(S^3)$ is a linear combination of elements built using at most two diagonals (and these are decomposable, in view of the codimension 5 case which we just treated), plus elements of the form
    \[ (p_{ij})^\ast(\Delta_S)\cdot (p_{k\ell})^\ast(\Delta_S)\cdot (p_{mn})^\ast(\Delta_S)\ .\]
    Since $(p_{12})^\ast(\Delta_S)\cdot (p_{13})^\ast(\Delta_S)=   (p_{12})^\ast(\Delta_S)\cdot (p_{23})^\ast(\Delta_S)$, these elements are of the form
    \[    D_{ij}:=  (p_{ij})^\ast(\Delta_S^2)\cdot (p_{k\ell})^\ast(\Delta_S)  \ .\]
    Since we know that $\Ss^{2/B_{dp}}\to B_{dp}$ has the Franchetta property (Proposition \ref{s2}), there is an equality of the form
      \[ \Delta_S^2 = \chi\, (p_1)^\ast(h^2)\cdot (p_2)^\ast(h^2)\ \ \ \hbox{in}\ A^4(S\times S) \]
      (where $\chi$ is the topological Euler characteristic of $S$). Using \eqref{excs}, we then find that
      \[  D_{ij}= \chi\,  (p_i)^\ast (h^2) \cdot (p_j)^\ast(h^2)\cdot (p_{k\ell})^\ast(\Delta_S)= \chi^\prime  (p_1)^\ast (h^2) \cdot (p_2)^\ast(h^2)  \cdot (p_3)^\ast(h^2)\ \ \ \hbox{in}\ A^6(S^3)\ ,\]
   for some $\chi^\prime\in\QQ$. It follows that all of $GDA^6(S^3)$ is decomposable, and hence injects into cohomology.
      \end{proof}

  \begin{proposition}\label{f3} Let $\YY\to B$ be as above. The family $\YY^{3/B}\to B$ has the Franchetta property in codimension $\le 2n$.
  \end{proposition}
  
  \begin{proof} The set-up $(\PP^{n+3}, \OO_{\PP^{n+3}}(2)^{\oplus 3})$ has the property $(\ast_3)$ of \cite{FLV}, and so the argument of loc. cit. implies that
    \begin{equation}\label{cube}  GDA^\ast_B(Y^3)= \langle h, \Delta_Y\rangle\ .\end{equation}
  
   In view of Proposition \ref{f2}, $GDA^j_B(Y^3)$ injects into cohomology for $j<2n$. We only need to treat the case of codimension $2n$. Note that \eqref{cube} shows in particular that
   the pullback $GDA^{2n}_{\MM^0_{(2,2,2)}}(Y^3)\to GDA^{2n}_B(Y^3)$ is an isomorphism. Thanks to Proposition \ref{key}, we are reduced to showing that $GDA^{2n}_{\MM_{dp}}(Y^3)$ injects into cohomology.
   
  For this, we exploit the relation between $Y$ and the associated double plane $S$ (Theorem \ref{net} and Corollary \ref{netchow}). The (generically defined) isomorphism of motives of Corollary \ref{netchow} induces a (generically defined) isomorphism
    \[ h(Y)\ \cong\ h(S)( 1-{n\over 2}) \oplus \bigoplus\one(\ast)\ \ \ \hbox{in}\ \MM_{\rm rat}\ .\]
    This induces a  
  commutative diagram 
    \[  \begin{array}[c]{ccc}  GDA^{2n}_{\MM_{dp}}(Y^3)& \xrightarrow{\cong}&  GDA^{3+{n\over 2}}_{\MM_{dp}}(S^3) \oplus  \bigoplus GDA^\ast_{\MM_{dp}}(S^2) \\
         &&\\
         \downarrow&&\downarrow\\
         &&\\
         H^{4n}(Y^3,\QQ) &\xrightarrow{\cong}& H^{6+n}(S^3,\QQ) \oplus \bigoplus H^\ast(S^2,\QQ) \\
         \end{array}\]
  
  In view of Proposition \ref{s2}, we only need the Franchetta property in codimension $3+{n\over 2}$ for $\Ss^{3/\MM_{dp}}\to \MM_{dp}$. 
  In case $n=2$, $\YY\to B$ is the universal family of genus 5 K3 surfaces; for this case, the proposition was already proven in \cite{FLV}. In case $n\ge 4$, the codimension $3+{n\over 2}$ is at least 5, and so the required Franchetta property for $S^3$
  follows from Proposition \ref{s3}.
    \end{proof}

  \section{Main result}
  
  \begin{theorem}\label{main} Let $Y\subset\PP^{n+3}$ be a smooth dimensionally transverse intersection of three quadrics, where $n\ge 2$ is even. Then $Y$ has a multiplicative Chow--K\"unneth decomposition.
  \end{theorem}

 \begin{proof} 
 
 Let us first construct a CK decomposition for $Y$. Letting $h\in A^1(Y)$ denote a hyperplane section, as before we consider the CK decomposition
    \[   \begin{split}  
                              \pi^{2j}_Y&:= {1\over 8}\, h^{n-j}\times h^j\ \ \ \ (j\not={n\over 2})\ ,\\
                                \pi^{n}_Y&:= \Delta_Y-\sum_{j\not={n\over 2}} \pi^{2j}_Y\ \ \ \ \ \ \in\ A^{n}(Y\times Y)\ .\\
                                \end{split}\]
 We observe that this CK decomposition is {\em generically defined\/} with respect to the family $\YY\to B$ (Notation \ref{not}), i.e. it is obtained by restriction from
 ``universal projectors'' $\pi^j_\YY\in A^{n-2}(\YY\times_B \YY)$. (This is just because $h$ and $\Delta_Y$ are generically defined.)    
 
  Writing $h^j(Y):=(Y,\pi^j_Y)\in\MM_{\rm rat}$, we have
   \[ h^{2j}(Y) \cong \one(-j)\ \ \  \hbox{in}\ \MM_{\rm rat}\ \ \ \ \ \ (j\not={n\over 2})\ ,\]
 i.e. the interesting part of the motive of $Y$ is contained in $h^{n}(Y)$.

 What we need to prove is that this CK decomposition is MCK, i.e.
      \begin{equation}\label{this} \pi_Y^k\circ \Delta_Y^{sm}\circ (\pi_Y^i\times \pi_Y^j)=0\ \ \ \hbox{in}\ A^{2n}(Y\times Y\times Y)\ \ \ \hbox{for\ all\ }i+j\not=k\ ,\end{equation}
      or equivalently that
       \[   h^i(Y)\otimes h^j(Y)\ \xrightarrow{\Delta_Y^{sm}}\  h(Y) \]
       coincides with 
       \[ h^i(Y)\otimes h^j(Y)\ \xrightarrow{\Delta_Y^{sm}}\ h(Y)\ \to\ h^{i+j}(Y)  \ \to\ h(Y)\ , \]   
       for all $i,j$.
       
     We observe that the cycles in \eqref{this} are generically defined, and that the vanishing \eqref{this} holds true modulo homological equivalence. As such, the required vanishing \eqref{this} follows from the Franchetta property for $\YY^{3/B}\to B$ in codimension $2n$ (Proposition \ref{f3}). 
       \end{proof}
 
 \begin{remark}\label{difficult} The argument proving Theorem \ref{main} fails for $n$ odd. The problem is that in case $n$ is odd, the primitive cohomology of $Y$ is related to a Prym variety \cite{Beau2} instead of the double plane of Theorem \ref{net}, and it seems difficult to establish the Franchetta property for these Prym varieties.

Also, the argument of Theorem \ref{main} fails for intersections $Y$ of four quadrics in an odd-dimensional projective space. In this case, the primitive cohomology of $Y$ is related to the primitive cohomology of a double solid \cite{Tera}. However, the double solids that appear in this construction are {\em not\/} general double solids of the given degree (this is due to the fact that
not every hypersurface in $\PP^3$ is determinantal), and so the methods of the present paper do not extend to this situation. 
 \end{remark}

  \section{Consequences}   
  
  \subsection{Chow ring of $Y$}
  
  \begin{corollary}\label{cor} Let $Y\subset\PP^{n+3}$ be as in Theorem \ref{main}. Then
  \[  \ima\Bigl( A^i(Y)\otimes A^j(Y)\to A^{i+j}(Y)\Bigr) =\QQ[h^{i+j}]=\QQ[c_{i+j}(Y)]\ \ \ \forall\ i,j>0\ .\]
\end{corollary}

\begin{proof} The equality $\QQ[h^{i+j}]=\QQ[c_{i+j}(Y)]$ is a formality: $Y\subset\PP^{n+3}$ is a complete intersection, and so (by adjunction) the Chow-theoretic Chern classes are proportional to powers of $h$.

The argument proving the equality $ \ima\bigl( A^i(Y)\otimes A^j(Y)\to A^{i+j}(Y)\bigr) =\QQ[h^{i+j}]$ is similar to \cite[Theorem 1.0.1]{Diaz}, which is about cubic hypersurfaces.

We begin by noting that when $n=2$, the complete intersection $Y$ is a K3 surface and so Corollary \ref{cor} is known from the seminal work of Beauville--Voisin \cite{BV}. We now suppose that $n\ge 4$, and so $Y$ is a Fano variety.

Let us consider the {\em modified small diagonal\/} 
  \[ \begin{split} \Gamma_3:= \Delta^{sm}_Y - {1\over 8}\Bigl( (p_{12})^\ast(\Delta_Y)\cdot (p_3)^\ast(h) +  (p_{13})^\ast(\Delta_Y)\cdot (p_2)^\ast(h)+  (p_{23})^\ast(\Delta_Y)\cdot (p_1)^\ast(h)\Bigr)&\\ +
    {\displaystyle\sum_{i+j+k=2n-4}} a_{ijk} (p_1)^\ast(h^i)\cdot (p_2)^\ast(h^j)\cdot (p_3)^\ast(h^k)\ \ \ \in\ A^{2n}(Y\times Y\times& Y)\ ,\\
    \end{split}\]
    which depends on choices of $a_{ijk}\in\QQ$. We make the following claim:
    
    \begin{claim}\label{cl} There exist $a_{ijk}\in\QQ$ such that
      \[  \Gamma_3=0\ \ \ \hbox{in}\ A^{2n}(Y\times Y\times Y)\]
      for all $Y$ as in Theorem \ref{main}.
      \end{claim}
      
      The claim implies the corollary, as can be seen by letting $\Gamma_3$ act on $\alpha\times\beta\in A^{i+j}(Y\times Y)$ (cf. \cite[Lemma 2.0.2]{Diaz}; this uses that $Y$ is Fano so that $A_0(Y)\cong\QQ$). To establish the claim, 
 we reason as in \cite[Proof of Proposition 2.8]{FLV3}: the MCK decomposition (plus the fact that $Y$ has transcendental cohomology, plus the relation of Lemma \ref{ok}) yields an identity of the form
      \[ \begin{split} \Delta^{sm}_Y={1\over 8}\Bigl( (p_{12})^\ast(\Delta_Y)\cdot (p_3)^\ast(h) +  (p_{13})^\ast(\Delta_Y)\cdot (p_2)^\ast(h)+  (p_{23})^\ast(\Delta_Y)\cdot (p_1)^\ast(h)\Bigr)&\\ + P\bigl( (p_1)^\ast(h),(p_2)^\ast(h), (p_3)^\ast(h)\bigr)\ \ \ \hbox{in}\ A^{2n}&(Y^3)\ ,\\
      \end{split}\]
   where $P$ is a symmetric polynomial with $\QQ$-coefficients.   
   
   (Alternatively, one could establish Claim \ref{cl} by using the Franchetta property for $Y^3$ in codimension $2n$ (Proposition \ref{f3}); thus one is reduced to finding $a_{ijk}$ such that
   the claim is verified modulo homological equivalence; this can be done as in \cite[Lemma 2.0.1]{Diaz}.)      
     \end{proof}

\begin{remark} Recall that when $n=2$, $Y$ is a K3 surface and Corollary \ref{cor} is immediate from the seminal work of Beauville--Voisin \cite{BV}.    

In case $n\ge 4$, the only non-trivial part of Corollary \ref{cor} is the fact that
   \begin{equation}\label{thisnot} \ima\Bigl( A^{n\over 2}(Y)\otimes A^1(Y)\to A^{{n\over 2}+1}(Y)\Bigr) = \QQ\ .\end{equation}
   (Indeed, it follows from Corollary \ref{netchow} or from \cite{Ot} that $A^j(Y)=\QQ[h^j]$ for all $j\not\in\{{n\over 2},{n\over 2}+1\}$; this covers all of Corollary \ref{cor} except for \eqref{thisnot}.)
 
 Let $Y^\prime\subset\PP^{2n+3}$ be the intersection of two quadrics in the net defining $Y$ (one may assume $Y^\prime$ is non-singular). Since $A^1(Y)\cong\QQ[h]$, one has an inclusion
   \begin{equation}\label{inc} \ima\Bigl( A^{n\over 2}(Y)\otimes A^1(Y)\to A^{{n\over 2}+1}(Y)\Bigr) \ \subset\ \ima\Bigl( A^{{n\over 2}+1}(Y^\prime)\to A^{{n\over 2}+1}(Y)\Bigr)\ ,\end{equation}
  and $A^{{n\over 2}+1}(Y^\prime)$ is infinite-dimensional. It seems reasonable to expect that \eqref{inc} is actually an equality (this would be the case if the inclusion morphism $Y\hookrightarrow Y^\prime$ is ``of pure grade 0'', in the sense of \cite{SV2}.)
   \end{remark}

 \subsection{The tautological ring}
 
 \begin{corollary}\label{cor1} Let $Y\subset\PP^{n+3}$ be a smooth complete intersection of three quadrics, where $n$ is even. Let $m\in\NN$. Let
  \[ R^\ast(Y^m):=\Bigl\langle (p_i)^\ast(h), (p_{ij})^\ast(\Delta_Y)\Bigr\rangle\ \subset\ \ \ A^\ast(Y^m)   \]
  be the $\QQ$-subalgebra generated by pullbacks of the polarization $h\in A^1(Y)$ and pullbacks of the diagonal $\Delta_Y\in A^{n}(Y\times Y)$. (Here $p_i$ and $p_{ij}$ denote the various projections from $Y^m$ to $Y$ resp. to $Y\times Y$).
  The cycle class map induces injections
   \[ R^\ast(Y^m)\ \hookrightarrow\ H^\ast(Y^m,\QQ)\ \ \ \hbox{for\ all\ }m\le 2b-1  \ ,\]
   where $b:=\dim H^{n}(Y,\QQ)$.
   
   Moreover, $R^\ast(Y^m)\to H^\ast(Y^m,\QQ)$ is injective for all $m$ if and only if $Y$ is Kimura finite-dimensional \cite{Kim}.
   \end{corollary}

\begin{proof} This is inspired by the analogous result for cubic hypersurfaces \cite[Section 2.3]{FLV3}, which in turn is inspired by analogous results for hyperelliptic curves \cite{Ta2}, \cite{Ta} (cf. Remark \ref{tava} below) and for K3 surfaces \cite{Yin}.

As in \cite[Section 2.3]{FLV3}, let us write $o:={1\over 8} h^{n}\in A^{n}(Y)$, and
  \[ \tau:= \Delta_Y - {1\over 8}\, \sum_{j=0}^{n}  h^j\times h^{n-j}\ \ \in\ A^{n}(Y\times Y) \]
  (this cycle $\tau$ is nothing but the projector on the motive $h^{n}_{prim}(Y)$).
Moreover, let us write 
  \[ \begin{split}   o_i&:= (p_i)^\ast(o)\ \ \in\ A^{n}(Y^m)\ ,\\
                        h_i&:=(p_i)^\ast(h)\ \ \in \ A^1(Y^m)\ ,\\
                         \tau_{ij}&:=(p_{ij})^\ast(\tau)\ \ \in\ A^{n}(Y^m)\ .\\
                         \end{split}\]
We now define the $\QQ$-subalgebra
  \[ \bar{R}^\ast(Y^m):=\Bigl\langle o_i, h_i, \tau_{ij}\Bigr\rangle\ \ \ \subset\ H^\ast(Y^m,\QQ)\ \ \ \ \ (1\le i\le m,\ 1\le i<j\le m)\ ; \]
 this is the image of $R^\ast(Y^m)$ in cohomology. One can prove (just as \cite[Lemma 2.11]{FLV3} and \cite[Lemma 2.3]{Yin}) that the $\QQ$-algebra $ \bar{R}^\ast(Y^m)$
  is isomorphic to the free graded $\QQ$-algebra generated by $o_i,h_i,\tau_{ij}$, modulo the following relations:
    \begin{equation}\label{E:X'}
			o_i\cdot o_i = 0, \quad h_i \cdot o_i = 0,  \quad 
			h_i^{n} =8\,o_i\,;
			\end{equation}
			\begin{equation}\label{E:X2'}
			\tau_{ij} \cdot o_i = 0 ,\quad \tau_{ij} \cdot h_i = 0, \quad \tau_{ij} \cdot \tau_{ij} = (b-1)\, o_i\cdot o_j
			\,;
			\end{equation}
			\begin{equation}\label{E:X3'}
			\tau_{ij} \cdot \tau_{ik} = \tau_{jk} \cdot o_i\,;
			\end{equation}
			\begin{equation}\label{E:X4'}
			\sum_{\sigma \in \mathfrak{S}_{b}} 
			\hbox{sign}(\sigma) 
			\prod_{i=1}^{b} \tau_{i, b+\sigma(i)} = 0\, ,
			\end{equation}

To prove Corollary \ref{cor1}, it suffices to check that these relations are verified modulo rational equivalence.
The relations \eqref{E:X'} take place in $R^\ast(Y)$ and so they follow from the Franchetta property for $Y$. 
The relations \eqref{E:X2'} take place in $R^\ast(Y^2)$. 
The first and the last relations are trivially verified modulo rational equivalence, because ($Y$ being Fano) one has
$A^{2n}(Y^2)=\QQ$. As for the second relation of \eqref{E:X2'}, this is the relation of Lemma \ref{ok}. 
   
   Relation \eqref{E:X3'} takes place in $\bar{R}^\ast(Y^3)$. Using the MCK decomposition (Theorem \ref{main}), we can verify that this relation is also verified modulo rational equivalence. Indeed, we have
   \[  \Delta_Y^{sm}\circ (\pi^{n}_Y\times\pi^{n}_Y)=   \pi^{2n}_Y\circ \Delta_Y^{sm}\circ (\pi^{n}_Y\times\pi^{n}_Y)  \ \ \ \hbox{in}\ A^{2n}(Y^3)\ ,\]
   which (using Lieberman's lemma) translates into
   \[ (\pi^{n}_Y\times \pi^{n}_Y\times\Delta_Y)_\ast    \Delta_Y^{sm}  =   ( \pi^{n}_Y\times \pi^{n}_Y\times\pi^{2n}_Y)_\ast \Delta_Y^{sm}                
                                            \ \ \ \hbox{in}\ A^{2n}(Y^3)\ ,\]
   which means that
   \[  \tau_{13}\cdot \tau_{23}= \tau_{12}\cdot o_3\ \ \ \hbox{in}\ A^{2n}(Y^3)\ .\]
   
  Finally, relation \eqref{E:X4'}, which takes place in $\bar{R}^\ast(Y^{2b})$
   is the Kimura finite-dimensionality relation \cite{Kim}: in the notation of loc. cit., relation \eqref{E:X4'} is written as   
    \[   \wedge^{b} \pi^{n, prim}_Y = 0\ \ \ \hbox{in}\ H^{b(2n)}(Y^{b}\times Y^{b},\QQ)\ ,\]
    which
  expresses the fact that $\dim H^{n}_{prim}(Y,\QQ)=b-1$ and so   
     \[ \wedge^{b} H^{n}_{prim}(Y,\QQ)=0\ .\]
 Assuming that $Y$ is Kimura finite-dimensional, this relation is also verified modulo rational equivalence.
   \end{proof}

\begin{remark}\label{tava} Given any curve $C$ and an integer $m\in\NN$, one can define the {\em tautological ring\/}
  \[ R^\ast(C^m):=  \Bigl\langle  (p_i)^\ast(K_C),(p_{ij})^\ast(\Delta_C)\Bigr\rangle\ \ \ \subset\ A^\ast(C^m) \]
  (where $p_i, p_{ij}$ denote the various projections from $C^m$ to $C$ resp. $C\times C$).
  Tavakol has proven \cite[Corollary 6.4]{Ta} that if $C$ is a hyperelliptic curve, the cycle class map induces injections
    \[  R^\ast(C^m)\ \hookrightarrow\ H^\ast(C^m,\QQ)\ \ \ \hbox{for\ all\ }m\in\NN\ .\]
   On the other hand, there are many (non hyperelliptic) curves for which the tautological rings $R^\ast(C^2)$ or $R^\ast(C^3)$ do {\em not\/} inject into cohomology (this is related to the non-vanishing of the Faber--Pandharipande cycle of Remark \ref{!}; it is also
   related to the non-vanishing of the Ceresa cycle, cf. \cite[Remark 4.2]{Ta} and \cite[Example 2.3 and Remark 2.4]{FLV2}). 
   
Corollary \ref{cor1} shows that the tautological ring of $Y$ behaves similarly to that of hyperelliptic curves, K3 surfaces and cubic hypersurfaces. 
\end{remark}

We can now improve on Proposition \ref{f3}, by removing the codimension $\le 2n$ restriction:   
   
\begin{proposition}\label{f3+} Let $\YY\to B$ be as in Notation \ref{not}. Then $\YY^{3/B}\to B$ has the Franchetta property.
\end{proposition}   
  
 \begin{proof} The set-up $(\PP^{n+3}, \OO_{\PP^{n+3}}(2)^{\oplus 3})$ has the property $(\ast_3)$ of \cite{FLV}. The argument of loc. cit. (in the precise form of \cite[Proposition 2.6]{FLV3}) then implies that
   \[ GDA^\ast_B(Y^3)=\langle h,\Delta_Y\rangle\ .\]
   The right-hand side injects into cohomology thanks to Corollary \ref{cor1}.
 \end{proof}

\subsection{The Chow ring of $S$}

\begin{theorem}\label{main2} Let $S$ be the double cover of $\PP^2$ branched along a smooth curve $C\subset\PP^2$. Then $S$ has an MCK decomposition.
 In particular,
   \[ \ima\Bigl( A^1(S)\otimes A^1(S)\to A^2(S)\Bigr) = \QQ[c_2(S)]\ .\]
\end{theorem}

\begin{proof} The degree of $C$ is necessarily even, and we may asume it is at least 6 (for otherwise, $S$ is Fano and there is nothing to prove). Let $h\in A^1(S)$ denote the pullback of a hyperplane section of $\PP^2$.
Setting
  \[   \begin{split}  \pi^0_S&:= {1\over 2}   \, h^2\times S\ ,\\
                             \pi^4_S&:= {1\over 2}  \, S\times h^2\ ,\\
                             \pi^2_S&:=\Delta_S-\pi^0_S-\pi^4_S\ \ \ \ \hbox{in}\ A^2(S\times S)\ ,\\
                             \end{split}\]
               one obtains a CK decomposition that is generically defined (with respect to $\MM_{dp}$).
             
 To verify that this CK decomposition is MCK, it suffices to know that $\Ss^{3/\MM_{dp}}\to \MM_{dp}$ has the Franchetta property. Via the bridge provided by Proposition \ref{key} and Corollary \ref{netchow}, this is implied by the Franchetta property for $\YY^{3/\MM^0_{(2,2,2)}}\to \MM^0_{(2,2,2)}$. Since $B\to \MM_{(2,2,2)}$ is dominant, this is in turn implied by the Franchetta property for $\YY^{3/B}\to B$ (Proposition \ref{f3+}).
             
 As for the second statement, this is a formal consequence of the MCK paradigm: indeed, all intersections of divisors are in $A^2_{(0)}(S)\cong\QQ$. The second Chern class $c_2(S)$ is generically defined (with respect to $B_{dp}$ or even $\MM_{dp}$), and so the Franchetta property for $\Ss\to B$ implies that $c_2(S)$ is also in $A^2_{(0)}(S)$.            
   \end{proof}

  \begin{remark} The K3-like behaviour of double planes expressed in Theorem \ref{main2} should be contrasted with the {\em non K3-like\/} behaviour of the surfaces exhibited by O'Grady \cite{OG1}: for any $r\in\NN$, there exists a surface $S\subset\PP^3$ (of degree very large) such that
  \[ \dim \ima \Bigl( A^1(S)\otimes A^1(S)\to A^2(S)\Bigr) \ \ge r\ .\]
  \end{remark}

  \begin{remark} To be honest, the last statement of Theorem \ref{main2} is not quite new. Indeed, the fact that 
    \begin{equation}\label{dim1}  \dim\, \ima\Bigl( A^1(S)\otimes A^1(S)\to A^2(S)\Bigr) =1 \end{equation} 
    for a double plane $S$
    can easily be seen directly: the Chow ring decomposes as a sum $A^\ast(S)=A^\ast(S)^+\oplus A^\ast(S)^-$ (these are the $+1$ resp. $-1$ eigenspaces with respect to the covering involution).
  One has $A^1(S)^+=\QQ[h]$, and (using that intersecting with $h\in A^1(S)$ is a multiple of $\tau^\ast \tau_\ast$, where $\tau\colon S\to \PP$ is the inclusion of $S$ into a weighted projective space $\PP$) it is readily checked that $A^1(S)^+\cdot A^1(S)^-=0$. Since $A^1(S)^-\cdot A^1(S)^-\subset A^2(S)^+=\QQ[h^2]$, the equality \eqref{dim1} follows.
  \end{remark}

\vskip1cm
\begin{nonumberingt} Thanks to "ik ben een aubergine-vliegje" Len.
\end{nonumberingt}

\end{document}